\documentclass[11pt]{amsart}

\usepackage{amsmath}
\usepackage{amssymb}

\newtheorem{theorem}{Theorem}[section]
\newtheorem{claim}[theorem]{Claim}
\newtheorem{mclaim}[theorem]{Main Claim}

\newtheorem{lemma}[theorem]{Lemma}

\newtheorem{construction}[theorem]{Construction}

\theoremstyle{definition}
\newtheorem{definition}[theorem]{Definition}

\theoremstyle{remark}

\newcount\skewfactor
\def\mathunderaccent#1#2 {\let\theaccent#1\skewfactor#2
\mathpalette\putaccentunder}
\def\putaccentunder#1#2{\oalign{$#1#2$\crcr\hidewidth
\vbox to.2ex{\hbox{$#1\skew\skewfactor\theaccent{}$}\vss}\hidewidth}}

\def\smallbox#1{\leavevmode\thinspace\hbox{\vrule\vtop{\vbox
   {\hrule\kern1pt\hbox{\vphantom{\tt/}\thinspace{\tt#1}\thinspace}}
   \kern1pt\hrule}\vrule}\thinspace}

\newcommand{\rest}{{\upharpoonright}}

\newcommand{\cf}{{\rm cf}}

\newcommand{\Dom}{{\rm Dom}}

\newcommand{\st}{{such that}}
\newcommand{\seq}{{sequence}}

\newcommand{\incr}{{increasing}}

\newcommand{\Wlog}{{Without loss of generality}}

\newcommand{\then}{{\underline{then}}}
\newcommand{\Then}{{\underline{Then}}}

\def\qedref#1{$\qed_{\reforiginal{#1}}$}

\title{Two cardinal models for singular $\mu$ }
\author{Shimon Garti}
\address{Institute of Mathematics The Hebrew University of 
Jerusalem Jerusalem 91904, Israel}
\email{shimonygarty@hotmail.com}
\thanks{The first author would like to thank Badri Kriheli, for the
encouragement.}

\author{Saharon Shelah}
\address{Institute of Mathematics
The Hebrew University of Jerusalem
Jerusalem 91904, Israel
and  Department of Mathematics
Rutgers University
New Brunswick, NJ 08854, USA}
\email{shelah@math.huji.ac.il}
\urladdr{http://shelah.logic.at}
\thanks{First typed: April 2006 \\
Research supported by the United States-Israel Binational
Science Foundation. Publication 891}

\begin{document}
\let\labeloriginal\label
\let\reforiginal\ref
\def\ref#1{\reforiginal{#1}}
\def\label#1{\labeloriginal{#1}}

\begin{abstract}
We deal here with colorings of the pair $(\mu^+,\mu)$,
when $\mu$ is a strong limit and singular cardinal. 
We show that there exists a coloring $c$, with no refinement.
It follows, that the properties of identities of 
$(\mu^+,\mu)$ when $\mu$ is singular, differ in 
an essential way from the case of regular $\mu$.
\end{abstract}

\maketitle

\newpage
 
\section{Introduction}

Identities (or identifications) were first defined by Shelah 
in the late 60-s. The purpose was dual. On the one 
hand, they may be used as a tool for solving problems in Model Theory. 
On the other hand, there is interest in them 
within the realm of Set theory.

The basic connection between identities and 
questions of model theory 
(especially the compactness question of various pairs of 
cardinals) or mathematical logic (like the subject of 
generalized quantifiers) is formulated in \cite{Sh:8}. 
It is used in a much more sophisticated context, in \cite{sh:604}.
But here, we are interested in pure set theoretical
considerations.

Shelah proved, in the first part of \cite{Sh:824} 
(i.e., \S0 and \S1), that the set of identities 
$ID_2(\mu^+,\mu)$ has the property of 2-simplicity.
He proved this, for a regular cardinal $\mu$, \st\ 
$\mu=\mu^{<\mu}$. A natural example is the pair 
$(\aleph_1,\aleph_0)$.

Now, one may ask if the assumption on $\mu$ is necessary.
We shall prove here, that it can hardly be avoided. 
We will take a singular $\mu$ \st\ $2^{<\mu}=\mu$. 
Even under that assumption, we will see that there exists 
$c:[\mu^+]^n\rightarrow \mu$ which is not 
computable from any coloring $d:[\mu^+]^m \rightarrow \mu$ 
when $m<n$.

Let us describe now the structure of the article. 
In section 1, we give some defintions and 
basic facts about identities. In section 2, 
we build the main proof under the assumption that
$2^\mu=\mu^+$. In section 3, we drop that 
assumtion, using methods of 
pcf theory. The result is that we have our theorem, 
even if $2^\mu>\mu^+$.

Let us try to explain the idea. Assume $\kappa=\cf(\mu)<\mu$.
Let $\langle \lambda_i:i<\kappa\rangle$ be an \incr\ \seq\ 
of regular cardinals, with limit $\mu$. 
Let $J=J^{\rm bd}_\kappa$ be the ideal of all the bounded subsets 
of $\kappa$. In section 3 we show that one 
can replace the assumption $2^\mu=\mu^+$ by the assumption that 
${\rm tcf}(\prod\limits_{i<\kappa}
\lambda_i,J)=\mu^+$. 

By Cohen, it is consistent with ZFC that 
$2^\mu\geq \aleph_\alpha$
for any ordinal $\alpha$, so $2^\mu>\mu^+$ is consistent, 
but not provable in ZFC. On the other hand, the 
fact that there exsists 
$\bar \lambda=\langle \lambda_i:i<\kappa\rangle$ \st\ 
${\rm tcf} (\prod\limits_{i<\kappa} \lambda_i,J)=\mu^+$ 
is a theorem of ZFC. 

That brings us to a philosophical question about the 
meaning of analizing the magnitude of $2^\mu$. It is 
clear that $2^\mu$ is easy to manipulate by forcing.
What do we do about this? In fact, several answers are 
possible. Pcf theory suggests that asking about the size of
$2^\mu$ is the wrong question. 

Instead of looking at the value of $2^\mu$, 
about which there is a vast variety of consistency results, 
we should ask the 
right questions about the cardinality of porducts
of cardinals, divided by an ideal. 
Section 3 here exemplifies the philosophical idea very well. 
Starting with a consistency result, we arrive at a real 
theorem of ZFC, by changing our focus from the continuum 
question to a statment about tcf. \newline
We thank deeply to the referee, for many helpful comments and improvements.
\newpage

\section{Definitions}

The basic notion that we need, is identity: 

\begin{definition}
\label{1.1}
\begin{enumerate}
\item[(a)] A partial identity ${\bf s}$ is a pair 
$(a,e)=(\Dom_{\bf s},e_{\bf s})$. 
$a$ is a finite set, and $e$ is an equivalence relation 
on a subfamily of the subsets of $a$.

We always require that $e$ respects the cardinality
of the subsets, i.e. $bec\Rightarrow |b|=|c|$.
\item[(b)] A full identity is an identity 
${\bf s}=(a,e)$, when $\Dom(e)={\mathcal P}(a)$.

We might say just ``identity'', instead of full 
identity.
\end{enumerate}

One may wonder, why do we distinguish between full
identities and partial identities. 
Well, in many cases we are interested in colorings 
of the type 
$c:[\lambda]^n\rightarrow \mu$ when $n$ is constant. 
Analizing those colorings helps us to understand identities 
with $e$ defined only on subsets of $a$ with cardinality $n$.
Those are partial identities, of course.
\end{definition}

\begin{definition}
\label{1.2}
Let $(a,e)$ be an identity (or a partial identity). 
We say that $\lambda\rightarrow (a,e)_\mu$ if 
for every function $f:[\lambda]^{<\aleph_0}\rightarrow \mu$
there is a one-to-one mapping $h:a\rightarrow \lambda$, \st\ 
$bec\Rightarrow f(h''(b))=f(h''(c))$. 

Notice, that the requirement of $\lambda\rightarrow 
(a,e)_\mu$ relates to every function $f$. 
So, the next definition which depends only 
on the pair $(\lambda,\mu)$, makes sense:
\end{definition}

\begin{definition}
\label{1.3}
$ID(\lambda,\mu):=\{(a,e):(a,e)$ is an identity, 
and $\lambda\rightarrow (a,e)_\mu\}$

But we might be interested also in the identities 
of a specific function $f$:
\end{definition}

\begin{definition}
\label{1.4}
Let $f:[\lambda]^{<\aleph_0}\rightarrow \mu$ be a function.\newline
$ID(f):=\{(a,e):(a,e)$ is an identity, and there exists a 
one-to-one mapping 
$$
h:a\rightarrow \lambda, \hbox{ \st\ } 
bec \Rightarrow f(h''(b))=f(h''(c))\}
$$
Notice that $ID(\lambda,\mu)\subseteq \bigcap \{ID(f):f$ 
is a function from $[\lambda]^{<\aleph_0}$ into $\mu\}$.\newline
One of the basic tools for investigating identities 
is the notion of refinement. 
The idea is to compute the values
of a coloring $c:[\lambda]^n\rightarrow \mu$, with 
a coloring $d:[\lambda]^m \rightarrow \mu$, when 
$m<n$.
\end{definition}

\begin{definition}
\label{1.5}
Let $m<n<\omega, (\lambda,\mu)$ a pair of infinite
cardinals. Let $c:[\lambda]^n \rightarrow \mu$ and 
$d:[\lambda]^m \rightarrow \mu$ be colorings. 
\begin{enumerate}
\item[(a)] We say that $d$ refines $c$, if:

For any $\alpha_0,\ldots,\alpha_{n-1}<\lambda$ with 
no repetitions, and any $\beta_0,\ldots,\beta_{n-1}<
\lambda$ with no repetitions, the condition
(*)\  is satisfied. This means \newline 
\begin{enumerate}
\item[(*)] If for every $u\in [n]^m$ we 
have $d(\{\alpha_\ell:\ell\in u\})=
d(\{\beta_\ell:\ell\in u\})$, \underline{then} 
$c(\{\alpha_0,\ldots,\alpha_{n-1}\})=
c(\{\beta_0,\ldots,\beta_{n-1}\})$.\newline
\end{enumerate}

\item[(b)] We say that $d$ is an order-
refinement for $c$, if we concentrate only
on the cases \st\ $\alpha_0<\ldots<\alpha_{n-1}$
and $\beta_0<\ldots<\beta_{n-1}$.
\end{enumerate}
\end{definition}

\newpage

\section{The main thoerem}

Let $\mu$ be a singular cardinal, 
$\mu=2^{<\mu}$. We deal, in this section, 
with the pair $(\mu^+,\mu)$. Through-out the whole
section we add the assumption that $2^\mu=\mu^+$. 

\begin{mclaim}
\label{2.1}
Assume: 
\begin{enumerate}
\item[(a)] $\mu$ is a singular cardinal
\item[(b)] $2^{<\mu}=\mu$
\item[(c)] $2^\mu = \mu^+$
\item[(d)] $n\in [2,\omega)$
\end{enumerate}
\Then\ there is a coloring $c:[\mu^+]^{n+1}
\rightarrow \mu$ \st\ no $d:[\mu^+]^n \rightarrow \mu$ 
is a refinement for $c$. Moreover, there is not 
even an order-refinement for $c$. 

{\rm Before beginning the proof, let us recall 
the parallel situation for a regular $\mu$. 
If $\mu=\mu^{<\mu}$, and 
$c:[\mu^+]^{<\aleph_0}\rightarrow \mu$ is a coloring, 
then there is $d:[\mu^+]^2\rightarrow \mu$ 
which is a refinement of $c$. We don't need the 
assumption of order on the ordinals in the 
domain of $c$. 

That theorem is the main claim in 
\cite[\S1]{Sh:824}. It follows, quite 
immediately, that $ID_2(\mu^+,\mu)$ is 
2-simple (Those notions are defined there).
So here we show that the colorings of 
$(\mu^+,\mu)$, when $\mu$ is singular, behave 
much differently. 

Let us go back to the claim. We shall start 
with a general lemma, which asserts the 
existence of a bounding function under 
some reasonable assumptions.} 
\end{mclaim}

\begin{lemma}
\label{2.2}
Let $\mu$ be a singular strong limit cardinal, 
and $n\in [2,\omega]$.

\Then\ we can find $\theta_n<\mu$ and 
$g_n:[\theta_n]^n\rightarrow \cf(\mu)$ 
\st:
\begin{enumerate}
\item[(*)] For every $f:[\theta_n]^{n-1}
\rightarrow \cf(\mu)$ there exists
$u_f\in [\theta_n]^n$ \st\ 
$v\in [u_f]^{n-1} \Rightarrow f(v)<g_n (u_f)$.
\end{enumerate}
\end{lemma}

\par \noindent Proof:\ 
Let $\kappa=\cf(\mu), \theta_2=\kappa^+$, 
and $\theta_{n+1}=\beth_{n-1} (\kappa^+)$ 
for every $n\in [2,\omega)$.

We separate the proof into two cases. 
In the first case $n=2$, and then we 
build directly the desired $g_2$, using 
the fact that $\kappa^+>\kappa$. In the 
second case we consider $n>2$, and we use 
an induction hypothesis.

\par \noindent 
\underline{Case 1}: 
$n=2$.

So we need $g_2:[\kappa^+]^2\rightarrow \kappa$, 
which dominates any $f:\kappa^+ \rightarrow \kappa$. 
For every $\alpha<\kappa^+$, let 
$h_\alpha:\alpha\rightarrow \kappa$ be a one-to-one 
mapping. Define for every $\alpha<\beta<\kappa^+ 
(=\theta_2)$ the following function: \newline
$$g_2 (\{\alpha,\beta\})=h_\beta (\alpha).$$

Let us try to show that $g_2$ is as required.
Assume that $f$ is a function from $\kappa^+$ 
into $\kappa$. By the pigeon hole principle, 
we can choose $\gamma<\kappa$ \st\ 
$S:=\{\alpha<\kappa^+:f(\alpha)=\gamma\}$ 
is of cardinality $\kappa^+$. 
We choose also an ordinal $\beta_*\in S$ 
\st\ $|S\cap \beta_*|=\kappa$. 

Notice that 
$$
|\{\{\alpha,\beta_*\}:
g_2 (\{\alpha,\beta_*\})\leq \gamma\}|=  
|\{\{\alpha,\beta_*\}:h_{\beta_*} 
(\alpha)\leq \gamma\}|<\kappa,
$$
since $\gamma<\kappa, \beta_*$ is constant 
and $h_{\beta_*}$ is one-to-one. But 
$|S\cap \beta_*|=\kappa$, so one may 
choose $\alpha_*\in S\cap \beta_*$ \st\ 
$g_2 (\{\alpha_*,\beta_*\})>\gamma$. 

On the other hand, $f(\alpha_*)=f(\beta_*)=\gamma$
(since both $\alpha_*$ and $\beta_*$ were
taken from $S$). Define 
$u_f=\{\alpha_*,\beta_*\}$, and we are done. 
\medskip

\par \noindent 
\underline{Case 2}: 
$n>2$. 

By the induction hypothesis, $\theta_\ell$ 
and $g_\ell:[\theta_\ell]^\ell\rightarrow \kappa$ 
satisfy the lemma for $\ell=n-1$.

Let $\langle f'_\alpha:\alpha\in 
[\theta_\ell,\theta_n)\rangle$ enumerate 
all the functions from $[\theta_n]^{n-1}$ 
into $\kappa$. Define $g_n:[\theta_n]^n
\rightarrow \kappa$ as follows. 
If $\alpha_0,\ldots, \alpha_{n-2}<\theta_\ell
\leq \alpha_{n-1}<\theta_n$, then let 
$g_n(\{\alpha_0,\ldots,\alpha_{n-1}\})$ be 
$$
{\rm max}\{f'_{\alpha_{n-1}} 
(\{\alpha_0,\ldots,\alpha_{n-2}\})+1,
g_\ell (\{\alpha_0,\ldots,\alpha_{n-2}\})\}.
$$
In any other case, let $g_n$ be zero. 

We will show that $(g_n,\theta_n)$ satisfies 
the claim. For this, assume $f$ is a function from 
$[\theta_n]^{n-1}$ into $\kappa$. Clearly, 
$f\rest [\theta_\ell]^{n-1}$ appears in the 
enumeration above. Let 
$\alpha_*\in [\theta_\ell,\theta_n)$ be an ordinal 
\st\ $f\rest [\theta_\ell]^{n-1}\equiv f'_{\alpha_*}$.
Define $f^- :[\theta_\ell]^{\ell-1} \rightarrow \kappa$ 
as follows:
$$
(\forall v\in [\theta_\ell]^{\ell-1})
(f^- (v)=f(v\cup \{\alpha_*\})).
$$
By the induction hypothesis, there exists 
$u_{f^-}=\{\alpha_0,\ldots,\alpha_{\ell-1}\}$ 
as required, i.e., if $v_m=u_{f^-}\setminus 
\{\alpha_m\}$ for every $m\leq \ell-1$ then 
$f^- (v_m)<g_\ell (u_{f^-})$. At last, we can define 
$u_f:= u_{f^-} \cup\{\alpha_*\}$. 
\begin{enumerate}
\item[$(*)_1$] $m\leq \ell-1 \Rightarrow 
f(v_m\cup \{\alpha_*\})= f^-(v_m)<g_\ell 
(u_{f^-})\leq g_n(u_{f^-} \cup \{\alpha_*\}) 
=g_n(u_f)$.\newline
\item[$(*)_2$] $f(u_{f^-})=f'_{\alpha_*} 
(u_{f^-})<g_n (u_{f^-} \cup \{\alpha_*\})=g_n(u_f)$.
\end{enumerate}
So, again, we are done. 
\hfill \qedref{2.2}
\medskip
${}$

Moving back to the the main claim, we try to create 
a coloring $c$ with no refinement. It is, somehow,
more convenient to work with functions that encode 
the information that the refinement captures, instead 
of dealing with the refinement itself. That's the idea 
behind the next lemma. 

\begin{lemma}
\label{2.3}
Let $\gamma\leq \mu^+$ be an ordinal, 
$c:[\gamma]^{n+1} \rightarrow \mu$ a coloring, and 
$d:[\gamma]^n\rightarrow \mu$ a refinement of $c$. 
\begin{enumerate}
\item[(a)] One can find $F:[\mu]^{n+1}\rightarrow \mu$
\st\ if $\alpha_0,\ldots,\alpha_n<\gamma$ with no 
repetitions, and for $0\leq \ell\leq n$ we write 
$d(\{\alpha_0,\ldots,\alpha_n\}\setminus
\{\alpha_\ell\})=\gamma_\ell<\mu$, \then\ 
$F(\{\gamma_0,\ldots,\gamma_n\})=c
(\{\alpha_0,\ldots,\alpha_n\})$.
\item[(b)] There exists $\gamma_*<\mu^+$ \st\ 
$F$ is definable from $d\rest [\gamma_*]^n$ 
and $c\rest [\gamma_*]^{n+1}$, even when 
$\gamma=\mu^+$.
\end{enumerate}
\end{lemma}

\par \noindent Proof:\ 
\begin{enumerate}
\item[(a)] Let $E$ be the equivalence relation that 
is determined by $c$, i.e. 
$$
\{\alpha_0,\ldots,\alpha_n\} E
\{\beta_0,\ldots,\beta_n\} \hbox{ iff }
c(\{\alpha_0,\ldots,\alpha_n\})=c
(\{\beta_0,\ldots,\beta_n\}).
$$
For any equivalence class of $E$,
choose a representative. If 
$\{\alpha_0,\ldots,\alpha_n\}\in [\gamma]^{n+1}$,
define $\gamma^\alpha_\ell=d(\{\alpha^*_0,
\ldots,\alpha^*_n\}\setminus \{\alpha^*_\ell\})$ when 
$\{\alpha^*_0,\ldots,\alpha^*_n\}$ is the representative
of the equivalence class $\{\alpha_0,\ldots,
\alpha_n\}/E$.\newline
Define $F(\{\gamma^\alpha_0,\ldots,\gamma^\alpha_n\})
=c(\{\alpha_0,\ldots,\alpha_n\})$ whenever 
$\{\alpha_0,\ldots,\alpha_n\}\in [\gamma]^{n+1}$. 
For every other $(n+1)$-tuple $\in[\mu]^{n+1}$, 
define $F$ to be zero. One can verify easily that 
$F$ is well defined and satisfies (a), because 
of the assumption that $d$ refines $c$. \newline
[Let us explain more thoroughly why $F$ is a well defined 
function from $[\mu]^{n+1}$ into $\mu$. Assume 
$\langle \gamma_0,\ldots,\gamma_n\rangle$ 
belongs to $[\mu]^{n+1}$. Choose a representative for 
every equivalence class of $E$. We split the definition 
into two cases. 

In the first case, there is no representative of 
the form $\{\alpha^*_0,\ldots,\alpha^*_n\}$ \st\ 
$$
d(\{\alpha^*_0,\ldots,\alpha^*_n\}\setminus \{\alpha^*_\ell\})
=\gamma_\ell \hbox{ for every } 0\leq \ell\leq n. 
$$ 
In that case we simply define 
$F(\{\gamma_0,\ldots,\gamma_n\})=0$.
Clearly, $F$ is well-defined in that case.

In the other case there is a representative 
$\{\alpha^*_0,\ldots,\alpha^*_n\}$ \st\ 
$$
d(\{\alpha^*_0,\ldots,\alpha^*_n\}\setminus 
\{\alpha^*_\ell\})=\gamma_\ell, 
\hbox{ for any } 0\leq \ell\leq n.
$$ 
We show that this representative is unique.
So suupose that $\{\beta^*_0,\ldots,\beta^*_n\}$ is 
also one of the $E$-representatives, and 
$d(\{\beta^*_0,\ldots,\beta^*_n\}\lambda 
\{\beta^*_\ell\})=\gamma_\ell$ for every 
$0\leq \ell\leq n$. It means that for every 
$u\in [n+1]^n$ we have 
$d(\{\alpha^*_\ell:\ell\in u\})=
d(\{\beta^*_\ell:\ell\in u\})$. 
By definition \ref{1.5}(a) we must infer that 
$c(\{\alpha^*_0,\ldots,\alpha^*_n\})=c
(\{\beta^*_0,\ldots,\beta^*_n\})$. By the definition 
of $E$ we have $\{\alpha^*_0,\ldots,\alpha^*_n\} E
\{\beta^*_0,\ldots, \beta^*_n\}$. But since we 
deal with representatives, and every equivalence class 
has only one representative, we conclude that 
$\alpha^*_\ell=\beta^*_\ell$ for every 
$0\leq\ell\leq n$. This fact enables us to define 
$F(\{\gamma_0,\ldots,\gamma_n\})=
c(\{\alpha^*_0,\ldots,\alpha^*_n\})$ 
(and even $=c(\{\alpha_0,\ldots,\alpha_n\})$), 
without any problem of ambiguity. So $F$ is well-defined also 
in that case].
\item[(b)] For every $\{\gamma_0,\ldots,
\gamma_n\}\in [\mu]^{n+1}$, choose 
$\{\alpha_0,\ldots,\alpha_n\}\in [\mu^+]^{n+1}$
\st\ $c(\{\alpha_0,\ldots,\alpha_n\})=
F(\{\gamma_0,\ldots,\gamma_n\})$, if there is such 
$\{\alpha_0,\ldots,\alpha_n\}$.
Define 
\[\begin{array}{ll}
\gamma_*:={\rm sup}\{\alpha_l+1:
(\exists\{\gamma_0,\ldots,\gamma_n\}\in 
[\mu]^{n+1}) (\exists\{\alpha_0,\ldots,
\alpha_n\}\in [\mu^+]^{n+1})\\ 
\qquad \qquad \ \{\alpha_0,\ldots,\alpha_n\}
\hbox{ was chosen as a witness}\\
\qquad \qquad\  \hbox{ for } 
\{\gamma_0,\ldots,\gamma_n\},\hbox{and}\ 0\leq l\leq n\}
\end{array}\]
Since $|[\mu]^{n+1}|=\mu$, and since 
$\mu^+$ is regular, we have 
$\gamma_*<\mu^+$. Clearly, 
$\gamma_*$ is as required in (b).
\end{enumerate}
\hfill \qedref{2.3}
\medskip
${}$

We are ready now to prove the theorem itself. 
First, we define a coloring $c$, by induction on 
$\alpha<\mu^+$. Then, we show that $c$ has no refinement. \newline

\begin{construction}
\label{2.4} 
${}$

{\rm Let $\langle \mu_\epsilon:\epsilon<\cf(\mu)\rangle$ 
be an \incr\ \seq\ of cardinals, with limit $\mu$. 
We define by induction on $\alpha<\mu^+,\alpha
\geq n+1$, the coloring $c_\alpha:
[\alpha]^{n+1}\rightarrow \mu$. We demand that 
$\beta<\alpha\Rightarrow c_\beta\subseteq c_\alpha$, 
so at the end we will be able to define 
$c=\bigcup\limits_{\alpha<\mu^+} c_\alpha$. 
Notice that $c\rest [\alpha]^{n+1}\equiv 
c_\alpha$, if we succeed. 

\par \noindent 
\underline{Stage 1}:
$\alpha=n+1$.

We define $c_\alpha (\{0,\ldots,n\})=0$
\medskip

\par \noindent 
\underline{Stage 2}: $\alpha=\beta+1$.

So $c_\beta$ was defined in the previous stage
and we need to build $c_{\beta+1}=c_\alpha$. 
Let $\Gamma=\{(d,\gamma):\theta_n\leq \gamma<\mu^+$,
and $d:[\gamma]^n\rightarrow \mu\}$, and let 
$\langle(d_\alpha,\gamma_\alpha):\alpha<\mu^+\rangle$ 
enumerate $\Gamma$ (remember that $2^\mu = \mu^+$).
\Wlog, if $d_\beta:[\gamma_\beta]^n\rightarrow \mu$
then $\gamma_\beta\leq \beta$. Also, let 
$\langle F_\alpha:\alpha<\mu^+,\alpha\geq \mu\rangle$
be an enumeration of the suitable functions $F$, 
obtained by virtue of Lemma \ref{2.3}. 

We try to find a pair $(d,\beta)$ \st\ $d$ 
refines $c_\beta$. If there is no refinement 
of $c_\beta$, we are done (extend $c_\beta$ trivially, 
and this gives a coloring with no refinement). 
So assume that there is $d:[\beta]^n\rightarrow \mu$ 
which refines $c_\beta$, and let $F:[\gamma_*]^{n+1}\rightarrow 
\mu$ be the function that computes $c_\beta$ out of 
the values of $d$ ($c_\beta$ exists, by Lemma \ref{2.3}(a) and (b)).

We would like to define $c_\alpha:[\alpha]^{n+1}
\rightarrow \mu$ by cases. 
Let $u\in [\alpha]^{n+1}$ be any $(n+1)$-tuple.
\begin{enumerate}
\item[(a)] If $u\subseteq \beta$, define 
$c_\alpha(u)=c_\beta(u)$.
\item[(b)] If $\beta\in u$, but $u\setminus 
\{\beta\}\nsubseteq \theta_n$, define 
$c_\alpha(u)=0$.
\item[(c)] If $\beta\in u$ and $u\setminus 
\{\beta\}\subseteq \theta_n$, define the 
set $W_u$ as follows:
$$
W_u:=\{F(\{\zeta_0,\ldots,\zeta_n\}):0
\leq \ell\leq n-1\Rightarrow \zeta_\ell<
\mu_{g_n (u\setminus \{\alpha\})} \hbox{ and }
\zeta_n=d(u\setminus \{\alpha\})\}
$$
clearly, $|W_u|<\mu$ (since 
$\mu_{g_n (u\setminus \{\alpha\})}<\mu$, 
and $\mu$ is a strong limit cardinal). 
Consequently, $\mu\setminus W_u\neq \emptyset$, so choose 
$\zeta\in \mu\setminus W_u$, and define $c_\alpha(u)=
\zeta$.
\end{enumerate}

\par \noindent 
\underline{Stage 3}: 
$\alpha$ is a limit ordinal.\newline
Define $c_\alpha=\bigcup\limits_{\beta<\alpha}c_\beta$.

Now, let $c=\bigcup\{c_\alpha:\alpha<\mu^+\}$. 
This gives a coloring $c:[\mu^+]^{n+1}\rightarrow 
\mu$.} 
\end{construction}
\hfill \qedref{2.4}

\begin{claim}
\label{2.5}
There is no $d:[\mu^+]^n\rightarrow \mu$ which refines $c$.
\end{claim}

\par \noindent Proof:\ 
Towards a contradiction, assume that $d$ refines $c$. 
We use $\gamma_*$ from Lemma \ref{2.3}, and $\theta_n$ 
from lemma \ref{2.2}. We may assume that $\theta_n
\leq \gamma_*$. 
We also use the enumeration of $\Gamma$, so we can 
find $\alpha<\mu^+$ \st\ $(d_\alpha,\gamma_\alpha) \equiv 
(d\rest [\gamma_*]^n, \gamma_*)$. 

Define $f:[\theta_n\setminus \{\alpha\}]^{n-1}
\rightarrow \cf(\mu)$ as follows:
$$
f(v)={\rm min}\{i:d(v\cup \{\alpha\})<\mu_i\}.
$$
As above, $F$ is the encoding of $d$. 
We may assume that $F\equiv F_\alpha$ in the 
enumeration of the $F$-s. By Lemma \ref{2.2}, 
there exists $u_f=\{\beta_0,\ldots,\beta_{n-1}\}$
\st\ $y\in [u_f]^{n-1}\Rightarrow f(y)<g_n (u_f)$. 

Take a closer look at $c(u_f\cup \{\alpha\})$. 
It was chosen in the $\alpha+1$-st stage, and 
its value is the value of 
$c_{\alpha+1} (u_f\cup \{\alpha\})$. 
We know that 
$$
c_{\alpha+1} (u_f\cup\{\alpha\})\neq
F_\alpha (d(u_f\setminus \{\beta_0\}\cup
\{\alpha\}),\ldots, d(u_f\setminus \{\beta_{n-1}\}
\cup \{\alpha\}),d(u_f))
$$
(provided that $d(u_f\setminus \{\beta_\ell\}
\cup \{\alpha\})<\mu_{g_n (u_f)}$ for $\ell<n$,
which holds here by Lemma \ref{2.2} and the definition of $f$).

In other words, $d$ fails to determine the value 
of $c$ on the set $u_f\cup \{\alpha\}$, contradicting the 
definition of $F$ from Lemma \ref{2.3}, which is based on the
assumption that $d$ refines $c$. 
\hfill \qedref{2.5},\qedref{2.1}
\medskip
\newpage

\section{the pcf advantage}

\begin{theorem}
\label{3.1}
Assume: 
\begin{enumerate}
\item[(a)] $\mu$ is a singular cardinal. 
\item[(b)] $2^{<\mu}=\mu$.
\item[(c)] $m\in [2,\omega)$.
\end{enumerate}
\Then\ there exists $c:[\mu^+]^{m+1}\rightarrow \mu$
\st\ no $d:[\mu^+]^m \rightarrow\mu$ refines it.
\end{theorem}

\par \noindent
Proof:\ 
Denote $\kappa=\cf(\mu)<\mu$, and 
$\theta=\theta_m=\beth_{m-2} (\kappa^+)$. 
Let $J=J^{\rm bd}_\kappa$, the ideal of bounded subsets of 
$\kappa$.  By \cite{Sh:g} (see Main Claim \ref{1.3} 
in Chapter II) we can choose an \incr\ \seq\ of regular 
cardinals $\langle \lambda_i:i<\kappa\rangle, 
\theta<\lambda_0$ and $\mu=\bigcup\limits_{i<\kappa} 
\lambda_i$, \st\ $\mu^+={\rm tcf} (\prod\limits_{i<\kappa}
\lambda_i,J)$.

Let $\langle g^*_\alpha:\alpha<\mu^+\rangle$ 
exemplify it. We may assume that the \seq\ of 
the $g^*_\alpha$-s is strictly \incr. We are 
going to define a coloring with no refinement, using 
the $g^*_\alpha$-s. But we need some other functions. 
\begin{enumerate}
\item[$(*)_0$] Let $\bar f^\theta=\langle f^\theta_\alpha
:\alpha<\mu^+\rangle$ be a \seq\ of functions \st:
\begin{enumerate}
\item[(a)] $f^\theta_\alpha:[\theta]^m\rightarrow \kappa$,
for any $\alpha<\mu^+$.
\item[(b)] For every $f:[\theta]^m\rightarrow\kappa$,
we have:
$$
\mu^+={\rm sup} \{\alpha:f^\theta_\alpha=f\}.
$$
\end{enumerate}
(The meaning of (b) is that every $f^\theta_\alpha$ 
appears $\mu^+$ times in the \seq. It enables us to 
pick a specific function from a high enough level in the 
\seq).
 \item[$(*)_1$] Let $h:[\theta]^m\rightarrow \kappa$ be 
a dominating function, as given in Lemma \ref{2.2}, i.e.,
for every $g:[\theta]^{m-1}\rightarrow \kappa$, 
there exists $v_g\in [\theta]^m$ \st:
$$
(\forall \gamma\in v_g) [g(v_g\setminus \{\gamma\})
<h(v_g)].
$$
Now, denote $n=m+1$, and define 
$c:[\mu^+]^n \rightarrow \mu$ as follows:
\begin{enumerate}
\item[(i)] If $v\in [\theta]^m$ and 
$\alpha\in [\theta,\mu^+)$, \then
$$
c(v\cup \{\alpha\}) := g^*_\alpha ({\rm max} 
\{h(v),f^\theta_\alpha (v)\}+1)+1.
$$
\item[(ii)] For any $u\in [\mu^+]^n$ that 
doesn't fall in (i), define $c(u)=0$.
\end{enumerate}
Assume towards a contradiction that 
$d:[\mu^+]^m\rightarrow \mu$ refines $c$. 
By Lemma \ref{2.3} we have $F:[\mu]^n\rightarrow \mu$
which computes $c$ from the values of $d$. 
We will reach the desired contradiction 
using $F$. 
We need some more functions:
\item[$(*)_2$] For every $j<\kappa$ and 
any $\alpha<\mu^+$, we define $f^*_{\alpha,j}:
[\theta]^{m-1}\rightarrow \kappa$ as follows:
$$
f^*_{\alpha,j} (v)={\rm Min}\{i<\kappa:i>j 
\hbox{ and } \lambda_i>d (v\cup \{\alpha\})\}.
$$
\item[$(*)_3$] Let $f^{**}:[\theta]^m\rightarrow \kappa$ 
be defined by:
$$
f^{**} (v)={\rm Min} \{i<\kappa:d(v)<\lambda_i\}.
$$
We add also some functions of a different form:
\item[$(*)_4$] Define $g'\in \prod\limits_{i<\kappa}
\lambda_i$ by 
$$
g'(i)={\rm sup} \{\lambda_i\cap {\rm Rang} 
(d\rest [\theta]^m)\}\cup \bigcup\limits_{j<i} \lambda_j.
$$
\item[$(*)_5$] Define $g''\in \prod\limits_{i<\kappa}
\lambda_i$ by 
$$
g''(i)=g'(i) \cup {\rm sup} \{\lambda_i\cap 
{\rm Rang} (F\rest [g'(i)]^n\}.
$$
\end{enumerate}
Everything is ready now. Since 
$g''\in \prod\limits_{i<\kappa} \lambda_i$, we 
can pick an ordinal $\alpha_0<\mu^+$ \st\ 
$g''<_J g^*_{\alpha_0}$. By $(*)_0$, we can choose 
$\alpha_0<\alpha_1<\mu^+$ \st\ 
$f^\theta_{\alpha_1}\equiv 
f^{**}$. 
Clearly, $g''<_J g^*_{\alpha_1}$, so by the nature of 
the ideal $J$, there exists $j(*)<\kappa$ \st\ 
$$
g'' \rest [j(*),\kappa)< g^*_{\alpha_1} \rest 
[j(*),\kappa).
$$
Choose $v_*\in [\theta]^m$ \st\ for every 
$\gamma\in v_*$ it is true that $f^*_{\alpha_1,j(*)} 
(v_*\setminus \{\gamma\})<h(v_*)$
(exists, by $(*)_1$). 
\relax From the definition of $f^*_{\alpha_1,j(*)}$, it follows that:
\begin{enumerate}
\item[$\odot_0$] $\gamma\in v_*\Rightarrow 
d (v_*\setminus \{\gamma\}\cup \{\alpha_1\})
<\lambda_{f^*_{\alpha_1,j(*)}(v_*\setminus
\{\gamma\})} <\lambda_{h(v_*)}$.

Let $i(*)={\rm max}\{h(v_*), f^\theta_{\alpha_1} (v_*)\}$. 
By the definition of the $f^*$-s, $\gamma\in v_*\Rightarrow
j(*)<f^*_{\alpha_1,j(*)} (v_*\setminus \{\gamma\})$, and 
since $f^*_{\alpha_1,j(*)} (v_*\setminus \{\gamma\})<h(v_*)$, 
we know that $j(*)<h(v_*)$. So $j(*)<i(*)$. 
We need this for bounding the values of the coloring $d$, 
because $\odot_0$ implies now that:
\item[$\odot_1$] $\gamma\in v_* \Rightarrow 
d (v_*\setminus \{\gamma\} \cup \{\alpha_1\})
<\lambda_{h(v_*)} \leq \lambda_{i(*)}$. 

This fact tells us what happens if we drop one ordinal
from $v_*$, adding $\alpha_1$ instead. 
We also know what happens 
if we omit $\alpha_1$ and keep $v_*$:
\item[$\odot_2$] $d (v_* \cup \{\alpha_1\} \setminus 
\{\alpha_1\})=d(v_*)< \lambda_{f^{**}(v_*)}=
\lambda_{f^\theta_{\alpha_1(v_*)}}\leq 
\lambda_{i(*)}$. \newline 
This follows by the definition of $f^{**}$ in $(*)_3$, and the nature 
of $\alpha_1$, which implies that $f^\theta_{\alpha_1}
\equiv f^{**}$.
\end{enumerate}
We can finish the proof now, as we did in the
former section. 
Define:
$$
W:= \{F(\zeta_0,\ldots,\zeta_{n-1}):\zeta_0<
\ldots<\zeta_{n-1}<\lambda_{i(*)}\},
$$
$$
W^+ := \{F(\zeta_0,\ldots,\zeta_{n-1}):
\zeta_0<\ldots<\zeta_{n-1}<g' (i(*)+1)\},
$$
and get $W\subseteq W^+$ and also $W\subseteq g'' (i(*)+1)$ 
(By $(*)_4$ and $(*)_5$).
By virtue of $F's$ definition we have $c(v_*\cup \{\alpha_1\})
\in g'' (i(*)+1)$.
On the other hand, 
$c(v_* \cup \{\alpha_1\})
=g^*_{\alpha_1} ({\rm max} \{h(v_*), f^\theta_{\alpha_1}
(v_*)\}+1)+1=g^*_{\alpha_1} (i(*)+1)+1$.
But $j(*)<i(*)+1$, so $g'' (i(*)+1)<g^*_{\alpha_1} 
(i(*)+1)+1$, a contradiction.
\hfill \qedref{3.1}
\medskip

Combine Theorem \ref{3.1} with the main claim of 
\cite[\S1]{Sh:824}, 
and one has (almost) a full picture for the 
pair $(\mu^+,\mu)$. 

We may wonder about the assumption $2^{<\mu}=\mu$. 
As a matter of fact, our proof depends only 
on the fact that $\theta=\beth_{m-2} (\kappa^+)<\mu$. 
Of course, we want this for every $m<\omega$,
but this is still a weaker assumption. 

We can also ask what happens for other pairs of 
cardinals. We will try, in a subsequent paper, 
to shed light on the pair $(\mu^{+n},\mu)$.

\end{document}